\documentclass[11pt]{amsart}
\usepackage{latexsym,amssymb,amsmath}
\usepackage{mathrsfs}
\usepackage{epsfig}
\usepackage{amsmath}
\usepackage{amsfonts}
\usepackage{amssymb}
\usepackage{amssymb,amscd}
\usepackage{graphicx}
\usepackage{subfigure}
\usepackage{psfrag}
\RequirePackage{graphics}

\begin{document}
\newtheorem{thm}{Theorem}
\newtheorem{lem}[thm]{Lemma}
\newtheorem{cor}[thm]{Corollary}
\newtheorem{conj}[thm]{Conjecture}
\newtheorem{qn}{Question}
\newtheorem{pro}{Proposition}[section]
\theoremstyle{definition}
\newtheorem{defn}{Definition}[section]
\newtheorem{remk}{Remark}[section]
\newcommand\bC{\mathbb{C}}
\newcommand\bR{\mathbb{R}}
\newcommand\bZ{\mathbb{Z}}
\def\square{\hfill${\vcenter{\vbox{\hrule height.4pt \hbox{\vrule width.4pt
height7pt \kern7pt \vrule width.4pt} \hrule height.4pt}}}$}

\newenvironment{pf}{{\it Proof.}\quad}{\square \vskip 12pt}

\title[PICARD MODULAR GROUP]{Generators of the Gauss-Picard modular groups in three complex dimensions }

\author{BaoHua Xie} \address{College of Mathematics and Econometrics \\Hunan University\\ Changsha, 410082, China}
\email{xiexbh@hnu.edu.cn}

\author{JieYan Wang} \address{College of Mathematics and Econometrics \\Hunan University\\ Changsha, 410082, China}
\email{jywang@hnu.edu.cn}

\author{YuePing Jiang} \address{College of Mathematics and Econometrics \\Hunan University\\ Changsha, 410082, China}
\email{ypjiang@hnu.edu.cn}



\begin{abstract}
In this paper, we prove that the  Gauss--Picard
modular group $\mathbf{PU}(3,1;\mathbb{Z}[i])$ in three complex dimensions can be generated by five
given transformations: two Heisenberg translations, two Heisenberg rotations and an involution. Indeed, our
method may work for the other higher dimensional Euclidean--Picard
modular groups.

\end{abstract}

\date{\today}

\maketitle
\section{Introduction}
The modular group $\mathbf{PSL}(2,\mathbb{Z})$  is a fundamental object of study in number theory, geometry, algebra, and many other areas of  mathematics.
There are some natural algebraic generalisation of the classical modular group. For a positive square-free integer $d$ , the Bianchi group $\Gamma_{d}$ is the group $\mathbf{PSL}(2, \mathcal{O}_d)$, where  $\mathcal{O}_d$ is the ring of integers in the imaginary quadratic number field  $\mathbb{Q}(\sqrt{-d})$.   A general method to determine finite presentation for each  $\Gamma_{d}$ was developed by Swan\cite{Sw} based on the geometrical work of Bianchi, while a separate purely algebraic method was given by Cohn \cite{Co}. As another  generalization of the modular group, the construction was generalised by Picard in 1883. Suppose that $\mathbf{H}$ is a Hermitian
matrix of signature $(2, 1)$ with entries in $\mathcal{O}_d$. Let $\mathbf{SU}(\mathbf{H};\mathcal{O}_d)$ denote the subgroup of $\mathbf{SU}(\mathbf{H})$ consisting of those matrices whose entries lie in $\mathcal{O}_d$.   Picard studied group $\mathbf{PU}(\mathbf{H};\mathcal{O}_d)$ acting on complex hyperbolic plane $\mathbf{H}_{\mathbb{C}}^{2}$. These groups, called Picard modular groups. Such groups have attracted a great deal of attention both for their intrinsic interest as discrete groups and also for their applications in complex hyperbolic geometry.

One can view  modular group  or Bianchi  group acting discontinuously on hyperbolic spaces.  Then Poincar\'{e}'s polyhedra Theorem provide a  geometric  method  to obtain their  generators  from their fundamental polyhedra.
But  Mostow's  work \cite{Mo}  told us that  the construction of a fundamental domain  in complex hyperbolic space  are rather more complicated than for spaces of constant curvature. Until recently, the geometry of  $\mathbf{SU}(\mathbf{H};\mathcal{O}_3)$ has been studied
by Falbel and Parker \cite{FP} and the geometry of $\mathbf{SU}(\mathbf{H};\mathcal{O}_1)$ has been studied  by
Francsics and Lax \cite{FL1} and Falbel, Francsics and Parker \cite{FFP}.
By applying similar idea of \cite{FP}, \cite{FFP}, Zhao\cite{Zhao} obtained generators of Euclidean--Picard groups
 $\mathbf{PU}(2,1;\mathcal{O}_d)$ for $d=2,7,11$.

There are also some simple algorithm allow us to obtain the generators of modular group or some Picard
modular group.  For example, the continued fraction algorithm may be applied to any element of the modular group $\mathbf{PSL}(2,\mathbb{Z})$.
This  shows that $S(z)=1/z$ and $T(z)=z+1$ generate $\mathbf{PSL}(2,\mathbb{Z})$.
This algorithm was extended to $\mathbf{PU}(2, 1;\mathbb{Z}[i])$ by Falbel et al.\cite{FFLP} thus giving a different system of
generators from that obtained via a fundamental domain in \cite{FFP}. In \cite{wxx}, the authors applied the continued fraction algorithm to $\mathbf{PU}(2, 1;\mathbb{Z}[\omega])$ and so produce a different generating system from that obtained in \cite{FP}.

There is an obvious
generalisation of Picard modular groups to higher complex dimensions. We note that very little is known about the geometry and algebraic properties, e.g. explicit fundamental domain, generators system of the higher dimensions Picard modular groups $\mathbf{PU}(n,1;\mathcal{O}_d)$.
In \cite{xwj},  the continued fraction algorithm had been generalised to Picard modular groups in higher complex dimensions.  It contained the first generalisation that we were aware of to a group of $4\times4$ matrices. However it seems very difficult to extend the continued fraction algorithm to
other higher dimensional Picard modular groups. In this paper, we obtain the generators of  Gauss--Picard
modular group $\mathbf{PU}(3,1;\mathbb{Z}[i])$ in three complex dimensions by using the geometric method of \cite{FP,FFP,Zhao}.

\section{Preliminaries}
\subsection{Complex Hyperbolic Space}
In this subsection, we recall some basic materials in complex hyperbolic geometry and Picard modular group.
The general reference on these topics are \cite{Go,Par}.

Let $\mathbb{C}^{n,1}$ denote the vector space $\mathbb{C}^{n+1}$ equipped with the Hermitian form
$$\langle \mathbf{w},\mathbf{z}\rangle=z_{1}\overline{w_{n+1}}+z_{2}\overline{w_{2}}+\ldots +z_{n}\overline{w_{n}}+z_{n+1}\overline{w_{1}}$$
where   $\mathbf{w}$ and $\mathbf{z}$ are the column vectors in $\mathbb{C}^{n,1}$ with entries $z_{1},z_{2},\ldots,z_{n},z_{n+1}$ and $w_{1},w_{2},\ldots,w_{n},w_{n+1}$ respectively. Equivalently, we may write
$$\langle \mathbf{w},\mathbf{z}\rangle=\mathbf{z}^{*}J\mathbf{w}$$ where $\mathbf{w}^{*}$ denote the Hermitian transpose of $\mathbf{w}$ and
\begin{equation*}
J=\left(\begin{array}{ccc}
0 & 0& 1\\
0 & I_{n-1}& 0\\
1 & 0& 0
\end{array}\right).
\end{equation*}

Consider the following subspaces of $\mathbb{C}^{n,1}$:
$$
V_{-}=\{\mathbf{v}\in \mathbb{C}^{n,1}: \langle \mathbf{v},\mathbf{v}\rangle<0\},
$$
$$
V_{0}=\{\mathbf{v}\in \mathbb{C}^{n,1}-\{0\}: \langle
\mathbf{v},\mathbf{v}\rangle=0\}.
$$

Let $\mathbb{P}: \mathbb{C}^{n,1}-\{0\}\rightarrow \mathbb{C}P^{n}$ be the canonical projection onto complex projective space.
Then the {\it complex hyperbolic $n$-space} is defined to be $\mathbf{H}_{\mathbb{C}}^{n}=\mathbb{P}(V_{-})$.  The boundary of the complex hyperbolic $n$-space $\mathbf{H}_{\mathbb{C}}^{n}$ consists of those points in $\mathbb{P}(V_{0})$ together with
a distinguished point at infinity, which denote $\infty$.
The finite points in the boundary of $\mathbf{H}_{\mathbb{C}}^{n}$ naturally
carry the structure of the generalized Heisenberg group (denoted by $\mathcal {H}_{2n-1}$),
which is defined to $\mathbb{C}^{n-1}\times \mathbb{R}$ with the group law
$$
(\xi,\nu)\cdot(z,u)=(\xi+z,\nu+u+2\Im\langle\langle \xi,z \rangle\rangle)
.$$
Here $\langle\langle \xi,z \rangle\rangle=z^{*}\xi$ is the standard positive defined Hermitian form on $\mathbb{C}^{n-1}$. In particular,
we write $\| \xi\|^{2}=\xi^{*}\xi$.

Motivated by this, we define horospherical coordinates on complex hyperbolic space.
To each point $(\xi,\nu, u)\in \mathcal {H}_{2n-1}\times \mathbb{R}_{+}$,
we associated a point $\psi(\xi,\nu, u)\in V_{-}$.
Similarly, $\infty$ and each point $(\xi,\nu, 0)\in \mathcal {H}_{2n-1}\times \{0\}$ is associated to a point in $V_{0}$ by $\psi$.
The map $\psi$ is given by
\begin{equation*}
\psi(\xi,\nu, u)=\left(\begin{array}{c}
(-|\xi|^{2}-u+i\nu)/2\\
\xi\\
1
\end{array}\right),\
\psi(\infty)=\left(\begin{array}{c}
1\\
0\\
\vdots\\
0
\end{array}\right).
\end{equation*}
We also define the origin $0$ to be the point in $\partial\mathbf{H}_{\mathbb{C}}^{n}$  with horospherical coordinates $(0,0,0)$.
We have
\begin{equation*}
\psi(0)=\left(\begin{array}{c}
0\\
0\\
\vdots\\
1
\end{array}\right).
\end{equation*}

The holomorphic isometry group of $\mathbf{H}_{\mathbb{C}}^{n}$ is the group $\mathbf{PU}(n,1)$ of complex
linear transformations, which preserve the above Hermintian form.
That is, for each element $G\in \mathbf{PU}(n,1)$, $G$ is unitary with respect to  $\langle \cdot,\cdot\rangle.$ The corresponding matrix $G=(g_{jk})^n_{i,j=1}$ satisfies the
following condition
\begin{equation}
G^*JG=J,
\end{equation}
where $G^*$ denote the conjugate transpose of the matrix $G$.
\subsection{Picard modular groups}
Let $\mathcal{O}_d$ be the ring of integers in the imaginary quadratic number field $\mathbb{Q}(i\sqrt{d})$ where
$d$ is a positive square free integer. If $d\equiv 1,2\mod4$, then $\mathcal{O}_d=\mathbb{Z}[\sqrt{d}i]$ and If $d\equiv 3\mod4$, then $\mathcal{O}_d=\mathbb{Z}[(1+\sqrt{d}i)/2]$. The subgroup of $\mathbf{PU}(n,1)$ with entries in $\mathcal{O}_d$ is called the Picard modular group for $\mathcal{O}_d$ and is written $\mathbf{PU}(n,1;\mathcal{O}_d)$. Obviously, if $d=1$, then the ring $\mathcal{O}_d$ can be written as $\mathbb{Z}[i]$.

 \begin{remk}
  The matrixes corresponding to the generators obtained in this paper belong to the group
$\mathbf{U}(3,1; \mathbb{Z}[i])$. In relation to complex hyperbolic isometries, the relevant group is $\mathbf{PU}(3,1; \mathbb{Z}[i])=\mathbf{SU}(3,1; \mathbb{Z}[i])/\mathbb{Z}/4$. The center of $\mathbf{SU}(3,1)$ is
isomorphic to $\mathbb{Z}/4$, the group of $4^{th}$ roots of unity.
 By abuse of notation, we will denote the Gauss--Picard modular group in three complex dimensions by  $\mathbf{U}(3,1; \mathbb{Z}[i])$.
\end{remk}

\subsection{Complex Hyperbolic Isometries}
   We now discuss the decomposition of complex hyperbolic isometries. We begin by considering those elements fixing $0$ and $\infty$.

The matrix group $\mathbf{U}(n-1)$ acts by Heisenberg rotation. In horospherical coordinates,
the action of $U\in\mathbf{U}(n-1)$ is given by
$$(\xi,\nu, u)\longmapsto (U\xi,\nu, u).$$
The corresponding matrix in $\mathbf{U}(n,1)$ acting on $\mathbb{C}^{n,1}$ is
\begin{equation*}
M_U\equiv\left(\begin{array}{ccc}
1 & 0& 0\\
0 & U& 0\\
0 & 0& 1
\end{array}\right).
\end{equation*}

The positive real numbers $r\in \mathbb{R}^{+}$ act by Heisenberg dilation. In horospherical coordinates, this acting is given by
$$(\xi,\nu, u)\longmapsto (r\xi,r^{2}\nu, r^{2}u).$$ In $\mathbf{U}(n,1)$ the corresponding matrix is
\begin{equation*}
A_r\equiv\left(\begin{array}{ccc}
r & 0& 0\\
0 & I_{n-1}& 0\\
0 & 0& 1/r
\end{array}\right).
\end{equation*}

The Heisenberg group acts by Heisenberg translation. For $(\tau,t)\in \mathcal {H}_{2n-1}$, this is
$$N_{(\tau,t)}(\xi,\nu)=(\tau+\xi, t+\nu+2\Im\langle\langle \tau, \xi \rangle\rangle).$$
As a matrix $N_{(\tau,t)}$ is given by
\begin{equation*}
N_{(\tau, t)}\equiv\left(\begin{array}{ccc}
1& -\tau^{*} & (-\|\tau\|^{2}+it)/2\\
0 & I_{n-1}& \tau\\
0 & 0& 1
\end{array}\right).
\end{equation*}

Heisenberg translations, rotations and dilations generate the Heisenberg similarity group. This is the full subgroup of $\mathbf{U}(n,1)$ fixing $\infty$.

Finally, there is one more important acting, called an inversion $R$, which interchanges $0$ and $\infty$. In matrix notation this map is
\begin{equation*}
R\equiv\left(\begin{array}{ccc}
0& 0& 1\\
0 & -I_{n-1}& 0\\
1 & 0& 0
\end{array}\right).
\end{equation*}

Let $\Gamma_{\infty}$ be the stabilizer subgroup of $\infty$ in $\mathbf{U}(n,1)$. That is
$$\Gamma_{\infty}\equiv\{g\in \mathbf{U}(n,1):\ g(\infty)=\infty \}.$$

\begin{lem}
Let $G=(g_{jk})^4_{j,k=1}\in \mathbf{U}(3,1)$. Then
$G\in\Gamma_{\infty}$ if and only if $g_{41}=0$.
\end{lem}

Using Langlands decomposition, any element
$P\in\Gamma_{\infty}$ can be decomposed as a product of a Heisenberg
translation, dilation, and a rotation:
\begin{equation}\label{DeofP}
P= N_{(\tau,t)} A_{r} M_{U} =\left(
            \begin{array}{ccc}
             r& -\tau^{*}U & (-\|\tau\|^2+it)/2r \\
              0 & U & \tau/r \\
              0 & 0 & 1/r \\
            \end{array}
          \right),
\end{equation}
The parameters satisfy the corresponding conditions. That is, $U\in \mathbf{U}(n-1),r\in \mathbb{R}^{+}$ and $(\tau,t)\in \mathcal {H}_{2n-1}$.

\medskip

\subsection{Isometric spheres}

\medskip

Given an element $G\in \mathbf{PU}(3,1)$ with satisfying $G(q_{\infty})\neq
q_{\infty}$, we define the isometric sphere of $G$ to be the
hypersurface
$$\left\{\mathbf{z}\in\textbf{H}^2_{\mathbb{C}}:|\langle \mathbf{z},q_{\infty}\rangle|=|\langle \mathbf{z},G^{-1}(q_{\infty})\rangle|\right\}.$$
For example, the isometric sphere of
$$I_0=\left[\begin{array}{cccc}
0&0&0&1\\
0&-1&0&0\\
0&0&-1&0\\
1&0&0&0
\end{array}\right]
$$
is \begin{equation}\label{eq:2-4}
\mathcal{B}_0=\left\{(\zeta_{1},\zeta_{2}, t,u)\in\mathfrak{S}:\left||\zeta_{1}|^2+|\zeta_{2}|^2+u+it\right|=2\right\}\end{equation} in horospherical
coordinates.

All other isometric spheres are images of $\mathcal
{B}_0$
by Heisenberg dilations, rotations and translations. Thus the
isometric sphere with radius $r$ and centre $(\zeta^{0}_{1},\zeta^{0}_{2},t^{0},0)$ is given
by
$$\left\{(\zeta_{1},\zeta_{2},t,u):\left||\zeta_{1}-\zeta^{0}_{1}|^2+|\zeta_{2}-\zeta^{0}_{2}|^2+u+it-it^{0}+2i\Im m(\zeta_{1}\bar{\zeta}_1^{0}+\zeta_{2}\bar{\zeta}_2^{0})\right|=r^2\right\}.$$

If $G$ has the matrix form \begin{equation}\label{eq:2-6}\left[\begin{array}{cccc}
a_{11}&a_{12}&a_{13}&a_{14}\\
a_{21}&a_{22}&a_{23}&a_{24}\\
a_{31}&a_{32}&a_{33}&a_{34}\\
a_{41}&a_{42}&a_{43}&a_{44}
\end{array}\right]
,\end{equation} then $G(q_{\infty})\neq q_{\infty}$ if and
only if $g\neq0$. The isometric sphere of $G$ has radius
$r=\sqrt{2/|a_{41}|}$ and centre $G^{-1}(q_{\infty})$, which in
horospherical coordinates is
$$(\zeta^{0}_{1},\zeta^{0}_{2},t^{0},0)=(\bar{a_{43}}/\bar{a_{41}},
\bar{a_{42}}/\bar{a_{41}},2\Im m(\bar{a_{44}}/\bar{a_{41}}),0).$$

\section{The generators of $\mathbf{U}(2; \mathcal{O}_1)$}

Let $\mathbf{U}(2; \mathcal{O}_1)$ be the unitary group $\mathbf{U}(2)$ over the ring $\mathcal{O}_1$. Recall that the
unitary matrix $A\in\mathbf{U}(2)$ is of the following form
\begin{equation*}
\mathbf{U}(2)=\left\{A=\left(\begin{array}{cc}
a & b\\
 -\lambda \overline{b}& \lambda \overline{a}
\end{array}\right):  |\lambda|=1, |a|^{2}+|b|^{2}=1 \right\}.
\end{equation*}
Then we can see that the elements in $\mathbf{U}(2; \mathcal{O}_1)$ are of the following form
\begin{equation*}
\left(\begin{array}{cc}
a & 0\\
 0& b
\end{array}\right),
\left(\begin{array}{cc}
0 & b\\
 a& 0
\end{array}\right)
\end{equation*} where $a,b$ are units in $\mathcal{O}_1$. Recall that the units of $\mathcal{O}_1$ are $\pm1,\pm i$.

 We can find that
\begin{equation*}
\left\{\left(\begin{array}{cc}
a & 0\\
 0& b
\end{array}\right): a,b=\pm 1,\pm i\right\}
\end{equation*}
 can be generated by
\begin{equation*}
\left(\begin{array}{cc}
1 & 0\\
 0&i
\end{array}\right),
\left(\begin{array}{cc}
i & 0\\
 0& 1
\end{array}\right).
\end{equation*} We also note that
\begin{equation*}
\left(\begin{array}{cc}
0 & 1\\
 1& 0
\end{array}\right)
\left(\begin{array}{cc}
a & 0\\
 0& b
\end{array}\right)=\left(\begin{array}{cc}
0 & b\\
 a& 0
\end{array}\right),
\end{equation*} and
\begin{equation*}
\left(\begin{array}{cc}
0 & 1\\
 1& 0
\end{array}\right)
\left(\begin{array}{cc}
i &0\\
0& 1
\end{array}\right)\left(\begin{array}{cc}
0 & 1\\
 1& 0
\end{array}\right)
=\left(\begin{array}{cc}
1& 0\\
 0& i
\end{array}\right).
\end{equation*}

Therefore we have the following result.
\begin{lem}
$\mathbf{U}(2; \mathcal{O}_1)$ can be generated by the following two unitary matrixes
\begin{equation*}
U_{1}=\left(\begin{array}{cc}
0 & 1\\
 1& 0
\end{array}\right),
 U_{2}=\left(\begin{array}{cc}
i &0\\
0& 1
\end{array}\right).
\end{equation*}
\end{lem}

\section{The  generators of the stabiliser}
Next, we consider the subgroup
$\Gamma_\infty$ of the Picard modular group
$\mathbf{U}(3,1;\mathcal{O}_1)$.

\begin{lem}Let
$\Gamma_{\infty}(3,1;\mathbb{Z}[i])$ denote the subgroup
$\Gamma_{\infty}$ of Picard modular group
$\mathbf{U}(3,1;\mathbb{Z}[i])$. Then any element
$P\in\Gamma_{\infty}(3,1;\mathbb{Z}[i])$ if and only if the
parameters in the Langlands decomposition of $P$ satisfy the
conditions
$$
r=1, t\in 2\mathbb{Z}, \tau=(\tau_{1},\tau_{2})^{T}\in\mathbb{Z}[i]^{2},
U \in \mathbf{U}(2; \mathbb{Z}[i])$$

\end{lem}
\begin{pf}
Let $P\in\Gamma_{\infty}(3,1;\mathbb{Z}[i])$ be the Langlands decomposition form (2). Then it is easy to see that $r=1$, $t\in 2\mathbb{Z}$ and $U \in \mathbf{U}(2; \mathbb{Z}[i])$.  Since the entries $\tau_{1},\tau_{2}$ of $\tau$ and the entry $(-\|\tau\|^2+it)/2$ are in the ring $\mathbb{Z}[i]$,
we get that   $|\tau_{1}|^2+|\tau_{2}|^2\in 2\mathbb{Z}$.
\end{pf}

\begin{pro} Let $\Gamma_{\infty}(3,1;\mathbb{Z}[i])$ be stated as above. Then $\Gamma_{\infty}(3,1;\mathbb{Z}[i])$  is generated
by the Heisenberg translation $N_{((1,1)^T,0)},N_{((0,0)^T,2)}$
and the Heisenberg rotations $M_{U_{i}}(i=1,2)$.
\end{pro}

\begin{pf}
 Suppose $P\in\Gamma_{\infty}(3,1;\mathbb{Z}[i])$. According to Lemma 3, there is no dilation
component in its Langlands decomposition, that is
$$
P=N_{(\tau,t)}M_U=\left(
  \begin{array}{ccc}
    1 & -{\tau}^* & (-||\tau||^2+it)/2 \\
    0 & I_{2} & \tau \\
    0 & 0 & 1 \\
  \end{array}
\right)\left(
         \begin{array}{ccc}
           1 & 0 & 0 \\
           0 & U& 0 \\
           0 & 0 & 1 \\
         \end{array}
       \right).
$$
Since the unitary matrix $U\in \mathbf{U}(2;\mathbb{Z}[i])$. Then the rotation component of $P$ in
the Langlands decomposition is generated by $M_{U_{i}}(i=1,2)$ by Lemma 3.

We now consider the Heisenberg translation part of $P$, $N_{(\tau,t)}$.
Let $$\tau=(m_{1}+n_{1}i, m_{2}+n_{2}i)^{T},$$ where $m_{1}, n_{1}, m_{2}, n_{2}\in\mathbb{Z}$, since
$\tau\in\mathbb{Z}[i]^{2}$. Since $|\tau|^{2}=m_{1}^{2}+n_{1}^{2}+ m_{2}^{2}+n_{2}^{2}\in 2\mathbb{Z},$
there are two case:

1). $m_{1}^{2}+n_{1}^{2}\in 2\mathbb{Z},m_{2}^{2}+n_{2}^{2}\in 2\mathbb{Z}$;

2). $m_{1}^{2}+n_{1}^{2}\in 2\mathbb{Z}+1,m_{2}^{2}+n_{2}^{2}\in 2\mathbb{Z}+1$.

We first consider the Case 1. We can write $\tau$ as follows $$\tau=\left(k_{1}(1+i)+l_{1}(1-i),k_{2}(1+i)+l_{2}(1-i)\right).$$
where $k_{1}, l_{1}, k_{2}, l_{2}\in\mathbb{Z}$. $N_{(\tau,t)}$ splits as
$$N_{(\tau,t)}=N_{((0,0)^{T},t)}\circ N_{(\tau,0)}.$$
Since $t=2k\in2\mathbb{Z}$, $N_{\left((0,0)^{T},t\right)}=N_{\left(0,0)^{T},2\right)}^{k}$.
We also have
$$N_{(\tau,0)}=N_{\left((1+i,0)^{T},0\right)}^{k_{1}}\circ N_{\left((i-1,0)^{T},0\right)}^{l_{1}} \circ N_{\left(0,0)^{T},2\right)}^{2k_{1}l_{1}}\circ
N_{\left((0,1+i)^{T},0\right)}^{k_{2}}\circ N_{\left((0,1+i)^{T},0\right)}^{l_{2}}\circ N_{\left(0,0)^{T},2\right)}^{-2k_{2}l_{2}}.$$

We also note  that
\begin{align*}
N_{\left((1+i,0)^{T},0\right)}&=N_{\left((1,1)^{T},0\right)}\circ N_{\left((i,-1)^{T},0\right)} \circ N_{\left(0,0)^{T},2\right)},\\
N_{\left((i-1,0)^{T},0\right)}&=N_{\left((i,1)^{T},0\right)}\circ N^{-1}_{\left((1,1)^{T},0\right)} \circ N^{-1}_{\left(0,0)^{T},2\right)},\\
N_{\left((0,1+i)^{T},0\right)}&=N_{\left((1,1)^{T},0\right)}\circ N_{\left((-1,i)^{T},0\right)} \circ N_{\left(0,0)^{T},2\right)},\\
N_{\left((0,i-1)^{T},0\right)}&=N_{\left((1,i)^{T},0\right)}\circ N^{-1}_{\left((1,1)^{T},0\right)}\circ N^{-1}_{\left(0,0)^{T},2\right)}.
\end{align*}
It is easy to see that
\begin{align*}
N_{\left((i,1)^{T},0\right)}&=M_{U_{2}}N_{\left((1,1)^{T},0\right)}   M^{-1}_{U_{2}} ,\\
N_{\left((i,-1)^{T},0\right)}&=M_{U_{1}} M^{2}_{U_{2}}M_{U_{1}}M_{U_{2}} N_{\left((1,1)^{T},0\right)}   M^{3}_{U_{2}} (M_{U_{1}} M_{U_{2}}M_{U_{1}})^{2},\\
N_{\left((-1,i)^{T},0\right)}&=M^{2}_{U_{2}} M_{U_{1}}M_{U_{2}}M_{U_{1}} N_{\left((1,1)^{T},0\right)}   M^{2}_{U_{2}} (M_{U_{1}} M_{U_{2}}M_{U_{1}})^{3},\\
N_{\left((1,i)^{T},0\right)}&=M_{U_{1}} M_{U_{2}}M_{U_{1}} N_{\left((1,1)^{T},0\right)}(M_{U_{1}} M_{U_{2}}M_{U_{1}})^{3}.
\end{align*}

In the Case 2, We only need to consider  the translation $N(\tau,0)\circ N_{\left((1,1)^{T},0\right)}$. So we will come back to Case 1.

\end{pf}

\section{Fundamental domain for the stabiliser in $\mathbf{PU}(2,1;\mathbb{Z}[i])$}
In \cite{FFP}, the authors described a method to find the fundamental domain for the
stabiliser of $q_{\infty}$ in Gauss-Picard modualr group $\mathbf{PU}(2,1;\mathbb{Z}[i])$ in two complex
dimensions. We review it now.

Let $\Gamma$ be $\mathbf{PU}(2,1;\mathbb{Z}[i])$ and $\Gamma_{\infty}$ be its stabiliser of $q_{\infty}$.
Every element of $\Gamma_\infty$ is upper triangular and its diagonal entries are units in $\mathbb{Z}[i]$. Recall that the units of $\mathbb{Z}[i]$ are $\pm1,\pm i$. Therefore $\Gamma_\infty$  contains no dilations and so is a subgroup of $Isom(\mathfrak{R})$ and fits into the exact sequence as
$$0\longrightarrow\mathbb{R}\cap\Gamma_\infty\longrightarrow
\Gamma_\infty\stackrel{\Pi_{*}}{\longrightarrow}\Pi_{*}(\Gamma_\infty)\longrightarrow1.$$

We  write the isometry group of the integer lattice as
\begin{equation*}
Isom(\mathbb{Z}[i])=\left\{\left[\begin{array}{cc}\alpha&\beta\\0&1\end{array}\right]:\alpha,\beta\in\mathbb{Z}[i], \alpha\ \text{is a unit}\right\}.
\end{equation*}

The following  proposition find the image and kernel in this exact sequence.
\begin{pro}\label{eq:3-1}The stabiliser
$(\Gamma_1)_\infty$ of $q_{\infty}$ in $\Gamma_1$ satisfies
$$0\longrightarrow \mathbb{Z}\longrightarrow(\Gamma_1)_\infty\stackrel{\Pi_{*}}{\longrightarrow}\Delta\longrightarrow1,$$
where $\Delta\subset Isom(\mathcal{O}_1)$ is of index 2.
\end{pro}

 As the first step toward the construction of a fundamental domain for the action of $\Gamma_{\infty}$ on $\mathfrak{R}$, one should  construct a fundamental domain in $\mathbb{C}$ of $\Delta\subseteq Isom(\mathbb{Z}[i])$.
 From the generators of $\Delta$, one can find that a fundamental domain for $\Delta\subseteq Isom(\mathcal{O}_1)$ is  the triangle $\vartriangle$ with vertices $0,1,i$; See Fig.1.

\begin{figure}
  \centering
  \includegraphics[]{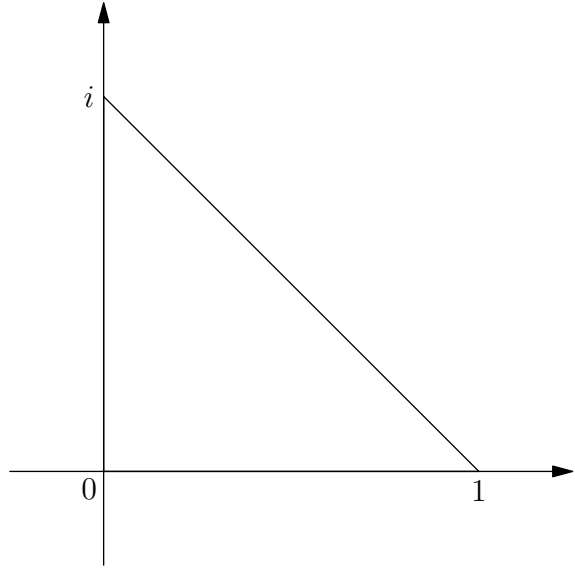}\\
  \caption{A Fundamental domain for the index two subgroup $\Delta \subset Isom(\mathbb{Z}[i])$ is a triangle $\vartriangle$ with vertices $0,1,i$.}
\end{figure}

In order to produce a fundamental domain for $\Gamma_{\infty}$ we look at all the preimages of
the triangle (that is a fundamental domain of $\Pi_*(\Gamma_{\infty})$)  under vertical projection
$\Pi$ and we intersect this with a fundamental domain for $ker(\Pi_*)$. The inverse of image of the
triangle under $\Pi$ is an infinite prism. The kernel of $\Pi_*$ is the infinite cyclic group
generated by $T$, the vertical translation by $(0,2)$. Hence a fundamental domain for $\Gamma_{\infty}$ is the
prism in $\mathfrak{R}$ with vertices $(0,\pm 1), (1,\pm 1),(i,\pm 1)$.

\section{The statement of our method and results}
In this section, we recall the geometric method used in (\cite{FP}, \cite{FFP}) to determine the generators of the Euclidean Picard groups and then state our method and results.

The geometric method is based on the special feature that the Euclidean-Picard modular orbifold has only one cusp for $d=1,2,7,11$. The basic idea of the proof can be described easily. Analogous to Theorem 3.5 of \cite{FP} one prove that $\langle \Gamma_\infty, R\rangle $ has only one cusp. The fact that $\mathbf{PU}(2,1;\mathcal{O}_d)$ has the same cusp and the stabiliser of infinity as the group generated by  $\langle \Gamma_\infty, R\rangle $  shows that they are the same. The key step is to find a union of isometric spheres such that a fundamental domain for $\Gamma_\infty$ is contained in the intersection of their exteriors and a fundamental domain for the stabiliser, which implies that the group  $\langle \Gamma_\infty, R\rangle $ has only one cusp. In other words, One should show that the union of the boundaries of these isometric spheres in Heisenberg group contains a fundamental domain for the stabiliser of infinity.

We will prove our theorems by using a similar idea with a litter different. The main observation is that there is no need to know the exact fundamental domain for the stabiliser of infinity. We will construct a set in heisenberg group which contain a fundamental domain for the stabiliser of infinity as a subset. Then we show that
the union of the boundaries of some isometric spheres in Heisenberg group covers this set. This also show that the group  $\langle \Gamma_\infty, R\rangle $ has only one cusp.

More prescisely, let $\Sigma$ be the following set
$$\{(\xi_{1},\xi_{2},t)|\xi_{i}\in\vartriangle, -1\leq t\leq 1 \}.$$ Here $\vartriangle$ is the fundamental domain of $\Delta \subset Isom(\mathbb{Z}[i])$.

Note that $\Sigma$ is not a fundamental domain for the stabiliser of infinity. Because some rotation preserves this set.
Indeed, the cone infinity from $\Sigma$ contains a fundamental domain for the stabiliser of infinity.
\begin{pro}The cone infinity from $\Sigma$ contains a fundamental domain for the stabiliser of infinity.
\end{pro}
\begin{pf}The restriction of the action of the stabiliser of infinity on each copy of $\mathbb{C}$ has the same fundamental doamin $\vartriangle$ of $\Delta\subset Isom(\mathbb{Z}[i])$. Then $\Sigma$ is the preimage
of $\vartriangle\times\vartriangle$ under vertical projection intersecting with a fundamental domain for
the vertical translation by $\left((0,0),2\right)$. Hence a fundamental domain for $\Gamma_{\infty}$ lies inside the cone infinity from $\Sigma$.
\end{pf}

In next section, we will prove the following theorem. Our main step is to show that $\Sigma$ lies inside the boundaries of some isometric spheres in Heisenberg group.

\begin{thm}  The Picard
modular group $\mathbf{U}(3,1;\mathbb{Z}[\omega_3])$ is generated by the Heisenberg translations
$$N_{\left((1,1)^{T},0\right)}=\left(\begin{array}{cccc}
             1& -1 & -1&-1\\
              0 & 1 & 0 &1\\
              0 & 0 & 1&1 \\
              0 & 0 & 0&1
            \end{array}
          \right), N_{\left((0,0)^{T},2\right)}=\left(\begin{array}{cccc}
             1& 0 & 0&i\\
              0 & 1 & 0 &0\\
              0 & 0 & 1&0 \\
              0 & 0 & 0&1
            \end{array}
          \right),$$ the Heisenberg rotations
         $$ M_{U_{1}}=\left(\begin{array}{cccc}
             1& 0 & 0&0\\
              0 &0 &1&0\\
              0 & 1&0&0 \\
              0 & 0 & 0&1
            \end{array}
          \right), 
M_{U_{2}}=\left(\begin{array}{cccc}
             1& 0 & 0&0\\
              0 & i & 0 &0\\
              0 & 0 & 1&0 \\
              0 & 0 & 0&1
            \end{array}
          \right),$$ and the  involution $$ R=\left(\begin{array}{cccc}
             0& 0 & 0&1\\
              0 & -1 & 0 &0\\
              0 & 0 &-1&0\\
              1& 0 & 0&0
            \end{array}
          \right).$$
\end{thm}

\section{The proof of the  theorem 4}
In this section, we will prove that the generators of Picard modular groups consist of the generators of the stabiliser and the involution.

Recall that the Cygan sphere $\mathcal{B}_{0}$ is the isometric sphere of $R$. The boundary of $\mathcal{B}_{0}$ is called the spinal sphere in Heisenberg group, we denote
by $S_{0}$ which  is defined by
$$\mathcal{S}_{0}=\{\left||\xi_{1}|^{2}+|\xi_{2}|^{2}+ti\right|
=2\}$$

Indeed we only need to consider the boundaries of isometric spheres in Heisenberg group because two isometric spheres have a non-empty interior
intersection if and only if the boundaries have a non-empty interior intersection.

It is not hard to see that party of $\Sigma$ lie outside $\mathcal{S}_{0}$. Therefore we need to find more isometric spheres whose
boundaries together with $\mathcal{S}_{0}$ contains the set $\Sigma$.

Note that $\Sigma$  has the following form
$$\Sigma=\{(\xi_{1},\xi_{2},t)|\xi_{1}\in\vartriangle,\xi_{2}\in\vartriangle, -1\leq t\leq 1\}.$$
First, we decompose $\vartriangle$ into three parts. That is
$\vartriangle=S_{1}\cup S_{2}\cup S_{3}$, where $S_{1}$ is a triangle with vertices $i,i/2, (1+i)/2$,
$S_{2}$ is a square with vertices $0,i/2, 1/2, (1+i)/2$,
$S_{3}$ is a triangle with vertices $0,1, (1+i)/2$. See Figure 2.

\begin{figure}
  \centering
  \includegraphics[]{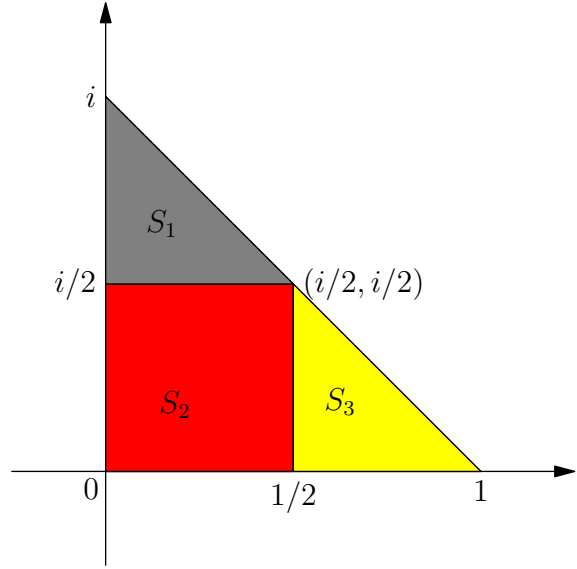}\\
  \caption{The decomposition of the fundamental domain $\vartriangle$ for $\Delta \subset Isom(\mathbb{Z}[i])$ into three parts.}
\end{figure}

Therefore $\Sigma$ will be decompose into nine subsets:

$\bullet$\  $\Sigma_{1}=\{(\xi_{1},\xi_{2},t)|\xi_{1}\in S_{1},\xi_{2}\in S_{1}, -1\leq t\leq 1\},$

$\bullet$\  $\Sigma_{2}=\{(\xi_{1},\xi_{2},t)|\xi_{1}\in S_{1},\xi_{2}\in S_{2}, -1\leq t\leq 1\},$

$\bullet$\  $\Sigma_{3}=\{(\xi_{1},\xi_{2},t)|\xi_{1}\in S_{1},\xi_{2}\in S_{3}, -1\leq t\leq 1\},$

$\bullet$\  $\Sigma_{4}=\{(\xi_{1},\xi_{2},t)|\xi_{1}\in S_{2},\xi_{2}\in S_{1}, -1\leq t\leq 1\},$

$\bullet$\  $\Sigma_{5}=\{(\xi_{1},\xi_{2},t)|\xi_{1}\in S_{2},\xi_{2}\in S_{2}, -1\leq t\leq 1\},$

$\bullet$\  $\Sigma_{6}=\{(\xi_{1},\xi_{2},t)|\xi_{1}\in S_{2},\xi_{2}\in S_{3}, -1\leq t\leq 1\},$

$\bullet$\  $\Sigma_{7}=\{(\xi_{1},\xi_{2},t)|\xi_{1}\in S_{3},\xi_{2}\in S_{1}, -1\leq t\leq 1\},$

$\bullet$\  $\Sigma_{8}=\{(\xi_{1},\xi_{2},t)|\xi_{1}\in S_{3},\xi_{2}\in S_{2}, -1\leq t\leq 1\},$

$\bullet$\  $\Sigma_{9}=\{(\xi_{1},\xi_{2},t)|\xi_{1}\in S_{3},\xi_{2}\in S_{3}, -1\leq t\leq 1\}.$

We first prove that $\mathcal{S}_{0}$ covers the subsets $\Sigma_{2}, \Sigma_{4}, \Sigma_{5}, \Sigma_{6}, \Sigma_{8}$.

If $(\xi_{1},\xi_{2},t)\in \Sigma_{5}$, then $|\xi_{1}|^{2}+|\xi_{2}|^{2}\leq (\frac{\sqrt{2}}{2})^{2}+(\frac{\sqrt{2}}{2})^{2}=1$.
So $\left||\xi_{1}|^{2}+|\xi_{2}|^{2}+it\right|\leq \sqrt{1+1}=\sqrt{2}<2$. Hence, $\Sigma_{5}\subset \mathcal{S}_{0}$.

If $(\xi_{1},\xi_{2},t)\in \Sigma_{2}$, $|\xi_{1}|^{2}+|\xi_{2}|^{2}\leq 1+\frac{\sqrt{2}}{2}=\frac{3}{2}$.
So $\left||\xi_{1}|^{2}+|\xi_{2}|^{2}+it\right|\leq \sqrt{(\frac{3}{2})^{2}+1}=\sqrt{\frac{13}{4}}<2$. Therefore, $\Sigma_{2}\subset \mathcal{S}_{0}$.

Similarly, we have $\Sigma_{4},\Sigma_{6},\Sigma_{8}$ are included in $\mathcal{S}_{0}$.

In order to prove this theorem, it is sufficient to prove that the remainian four subsets are covered by some Heisenberg spheres.

For the set $\Sigma_{9}$, we consider the map $N_{\left((1,1)^{T},0\right)}RN^{-1}_{\left((1,1)^{T},0\right)}$ whose isometric sphere which we denote by $\mathcal{B}_{1}$ is Cygan sphere centred at
the point $\left((1,1)^{T},0,0\right)$(in horospherical coordinates) with radius 1. The boundary of $\mathcal{B}_{1}$ is a Heisenberg sphere given by
$$\mathcal{S}_{1}=\{\left||\xi_{1}-1|^{2}+|\xi_{2}-1|^{2}+i(t+2\Im(\xi_{1}+\xi_{2}))\right|
=2\}$$
If $(\xi_{1},\xi_{2},t)\in \Sigma_{9}$, then $\xi_{1}\in S_{3},\xi_{2}\in S_{3},-1\leq t \leq 1$. We get that
$$0\leq \Im \xi_{i}\leq \frac{1}{2}, |\xi_{i}-1|^{2}\leq \frac{1}{2}.$$ So $$-1\leq t+2\Im(\xi_{1}+\xi_{2})\leq 3.$$
It is easy to see that the subset
$$\{(\xi_{1},\xi_{2},t)|\xi_{1}\in S_{3},\xi_{2}\in S_{3}, -1\leq t+2\Im(\xi_{1}+\xi_{2})\leq 1\}$$
lies inside $\mathcal{S}_{1}$ and the set $$\{(\xi_{1},\xi_{2},t)|\xi_{1}\in S_{3},\xi_{2}\in S_{3}, 1\leq t+2\Im(\xi_{1}+\xi_{2})\leq 3\}$$ lies inside $$T^{-1}(\mathcal{S}_{1})=\{\left||\xi_{1}-1|^{2}+|\xi_{2}-1|^{2}+i(t-2+2\Im(\xi_{1}+\xi_{2}))\right|
=2\}.$$
Therefore, $S_{1}$  and $T^{-1}(S_{1})$ cover the set $\Sigma_{9}$.

For the set $\Sigma_{7}$, we consider the map $N_{\left((1,i)^{T},0\right)}RN^{-1}_{\left((1,i)^{T},0\right)}$ whose isometric sphere which we denote by $\mathcal{B}_{2}$ is Cygan sphere centred at
the point $\left((1,i)^{T},0,0\right)$. The boundary of $\mathcal{B}_{2}$ is given by
$$\mathcal{S}_{2}=\{\left||\xi_{1}-1|^{2}+|\xi_{2}-i|^{2}+i(t+2\Im(\xi_{1})+2\Re(\xi_{2}))\right|
=2\}.$$

If $(\xi_{1},\xi_{2},t)\in \Sigma_{7}$, then $\xi_{1}\in S_{3},\xi_{2}\in S_{1},-1\leq t \leq 1$. We get that
$$0\leq \Im \xi_{1}\leq \frac{1}{2},0\leq \Re \xi_{2}\leq \frac{1}{2}, |\xi_{1}-1|^{2}\leq \frac{1}{2},|\xi_{2}-i|^{2}\leq \frac{1}{2}.$$ So $$-2\leq t+2\Im(\xi_{1}+\xi_{2})\leq 2.$$

If $-1\leq t+2\Im(\xi_{1}+\xi_{2})\leq 1$, then the subset
$$\{(\xi_{1},\xi_{2},t)|\xi_{1}\in S_{3},\xi_{2}\in S_{1}, -1\leq t+2\Im(\xi_{1}+\xi_{2})\leq 1\}$$ lies inside $\mathcal{S}_{2}$.

If $-2\leq t+2\Im(\xi_{1}+\xi_{2})\leq -1$, then the subset
$$\{(\xi_{1},\xi_{2},t)|\xi_{1}\in S_{3},\xi_{2}\in S_{1}, -2\leq t+2\Im(\xi_{1}+\xi_{2})\leq -1\}$$ lies inside $T(\mathcal{S}_{2})$.

If $1\leq t+2\Im(\xi_{1}+\xi_{2})\leq 2$, then the subset
$$\{(\xi_{1},\xi_{2},t)|\xi_{1}\in S_{3},\xi_{2}\in S_{1}, 1\leq t+2\Im(\xi_{1}+\xi_{2})\leq 2\}$$ lies inside $T^{-1}(\mathcal{S}_{2})$.

For the set $\Sigma_{1}$, we consider the map $N_{\left((i,i)^{T},0\right)}RN^{-1}_{\left((i,i)^{T},0\right)}$ whose isometric sphere which we denote by $\mathcal{B}_{3}$ is Cygan sphere centred at
the point $\left((i,i)^{T},0,0\right)$. The boundary of $\mathcal{B}_{3}$ is a Heisenberg sphere given by
$$\mathcal{S}_{3}=\{\left||\xi_{1}-1|^{2}+|\xi_{2}-1|^{2}+i(t+2\Im(\xi_{1}+\xi_{2}))\right|
=2\}.$$

If $(\xi_{1},\xi_{2},t)\in \Sigma_{1}$, then $\xi_{1}\in S_{1},\xi_{2}\in S_{1},-1\leq t \leq 1$. We get that
$$0\leq \Re\xi_{i}\leq \frac{1}{2}, |\xi_{i}-1|^{2}\leq \frac{1}{2}.$$ So $$-3\leq t+2\Im(\xi_{1}+\xi_{2})\leq 1.$$
As before, we can see that
$\Sigma_{1}$ covered by the Heisenberg spheres corresponding to the maps $N_{\left((i,i)^{T},0\right)}RN^{-1}_{\left((i,i)^{T},0\right)}$ and $TN_{\left((i,i)^{T},0\right)}RN^{-1}_{\left((i,i)^{T},0\right)}T^{-1}$.

Lastly, we consider the set  $\Sigma_{3}$, we consider the map $N_{\left((i,1)^{T},0\right)}RN^{-1}_{\left((i,1)^{T},0\right)}$ whose isometric sphere
which we denote by $\mathcal{B}_{4}$ is Cygan sphere centred at
the point $\left((i,1)^{T},0,0\right)$. The boundary of $\mathcal{B}_{4}$ is a Heisenberg sphere given by
$$\mathcal{S}_{4}=\{\left||\xi_{1}-i|^{2}+|\xi_{2}-1|^{2}+i(t-2\Re(\xi_{1})+2\Im(\xi_{2}))\right|
=2\}.$$

If $(\xi_{1},\xi_{2},t)\in \Sigma_{1}$, then $\xi_{1}\in S_{1},\xi_{2}\in S_{3},-1\leq t \leq 1$. We get that
$$0\leq \Re\xi_{1}\leq \frac{1}{2},0\leq \Im\xi_{1}\leq \frac{1}{2}, |\xi_{1}-i|^{2}\leq \frac{1}{2},|\xi_{2}-1|^{2}\leq \frac{1}{2}.$$ So $$-2\leq t-2\Re(\xi_{1})+2\Im(\xi_{2})\leq 2.$$

If $-1\leq t-2\Re(\xi_{1})+2\Im(\xi_{2})\leq 1$, then the subset
$$\{(\xi_{1},\xi_{2},t)|\xi_{1}\in S_{3},\xi_{2}\in S_{1}, -1\leq t-2\Re(\xi_{1})+2\Im(\xi_{2})\leq 1\}$$ lies inside $\mathcal{S}_{4}$.

If $-2\leq t-2\Re(\xi_{1})+2\Im(\xi_{2})\leq -1$, then the subset
$$\{(\xi_{1},\xi_{2},t)|\xi_{1}\in S_{3},\xi_{2}\in S_{1}, -2\leq t-2\Re(\xi_{1})+2\Im(\xi_{2})\leq -1\}$$ lies inside $T(\mathcal{S}_{4})$.

If $1\leq t-2\Re(\xi_{1})+2\Im(\xi_{2})\leq 2$, then the subset
$$\{(\xi_{1},\xi_{2},t)|\xi_{1}\in S_{3},\xi_{2}\in S_{1}, 1\leq t-2\Re(\xi_{1})+2\Im(\xi_{2})\leq 2\}$$ lies inside $T^{-1}(\mathcal{S}_{4})$.
Thus, $\Sigma_{3}$ covered by the Heisenberg spheres corresponding to the maps $$N_{\left((i,1)^{T},0\right)}RN^{-1}_{\left((i,1)^{T},0\right)}, TN_{\left((i,1)^{T},0\right)}RN^{-1}_{\left((i,1)^{T},0\right)}T^{-1}, T^{-1}N_{\left((i,1)^{T},0\right)}RN^{-1}_{\left((i,1)^{T},0\right)}T.$$

\begin{remk}Our method may work for the other higher dimensional Euclidean--Picard
modular groups. But the calculation will be more complicated. For example,
 The set $\Sigma$ will be decomposed into smaller parts. Then one need more Heisenberg spheres
 to cover the set $\Sigma$.
\end{remk}

\section*{Acknowledgements}

This work was partially supported by NSF(No.11071059).
B. Xie also supported by NSF(No.11201134) and 'Young teachers support program' of Hunan University.

\end{document}